\documentstyle [12pt]{article}

\title {On special partial types and weak canonical bases in simple theories}
\author {Ziv Shami}
\newtheorem {theorem}{Theorem}[section]
\newtheorem {lemma}[theorem]{Lemma}
\newtheorem {definition}[theorem]{Definition}

\newtheorem {fact}[theorem]{Fact}

\newtheorem {corollary}[theorem]{Corollary}

\newtheorem {remark}[theorem]{Remark}

\newtheorem {claim}[theorem]{Claim}

\newtheorem {example}[theorem]{Example}

\def\proof {\noindent \textbf{Proof:} }

\newsavebox{\indbin}
\savebox{\indbin}{\begin{picture}(0,0)
\newlength{\gnu}
\settowidth{\gnu}{$\smile$} \setlength{\unitlength}{.5\gnu} \put(-1,-.65){$\smile$}
\put(-.25,.1){$|$}
\end{picture}}
\newcommand{\nonfork}[3]
{\mbox{$\begin{array}{ccc} \mbox{$#1$} & \usebox{\indbin} & \mbox{$#2$} \\
        & \mbox{$#3$} &
\end{array}$}}
\newcommand{\nonforkempty}[2]
{\mbox{$\begin{array}{ccc} \mbox{$#1$} & \usebox{\indbin} & \mbox{$#2$}
\end{array}$}}

\newcommand{\forkempty}[2]
{\mbox{$\begin{array}{ccc} \mbox{$#1$} & \!\mbox{$\!\!\not\!\:\usebox{\indbin}$} & \mbox{$#2$}
\end{array}$}}

\def\card #1 {{\vert #1 \vert}}

\def\CC {{\cal C}}

\def\UU {{\cal U}}

\begin{document}
\maketitle


\section{Introduction}
For a simple theory the notion of the canonical base is essential for the development of parts of
the theory such as the theory of analyzability. Given an amalgamation base $p\in S(A)$, the
canonical base of $p$ is the minimal hyperimaginary, in the sense of definable closure, $e\in
dcl(A)$ such that $p$ doesn't fork over $e$ and $p\vert e$ is an amalgamation base. In this note we
define a notion of a weak canonical base for a partial type in a simple theory; it is defined in
the same way as the usual canonical base except that it is required to be minimal with respect to
bounded closure in the above sense (and there is no requirement on the restriction of the partial
type to it). We prove that members of a certain family of partial types (we call them special
partial types) have a weak canonical base. This family clearly properly contains the class of
amalgamation bases. Our original motivation was to prove Corollary \ref{cor1} for obtaining certain
definability result that seemed required for the proof of the dichotomy between 1-basedness and
supersimplicity proved in [S1]; however, this corollary turned out to be unnecessary for this
specific definability result. Nevertheless, it should have other applications to situations where
one needs a compactness argument when dealing with certain family of canonical bases.

The characterization of the class of partial types that admit a weak canonical base appears to be
an important problem and it looks reasonable that this class should properly contain the class of
special partial types that we deal with in this paper. The class of special partial types is a
certain subclass of the class of partial types obtained by generic composition of a pair of
complete types; we say that a partial type $r(x,a)$ over a sufficiently saturated model of $T$ is
obtained by generic composition of the complete types $p(x,y)$ and $q(y,z)$ (without parameters) if
the following condition holds in that model: $b$ realizes $r$ iff there exists $c$ such that
$p(b,c)$, $q(c,a)$, and $\nonfork{b}{a}{c}$. Our proof doesn't seem to extend for general generic
composition. The skeleton of the proof of the existence of weak canonical bases is similar to the
construction of the usual canonical base. Throughout this paper, $T$ is assumed to be a first-order
simple theory and we work in a monster model $\CC$ of $T$, namely a sufficiently saturated, and
sufficiently strongly-homogeneous model of $T$. We will sometime assume, for simplicity, that $T$
is hypersimple, namely, a simple theory with elimination of hyperimaginaries. We only assume basic
knowledge of simple theories as in [K],[KP] and [HKP].

\section{Weak canonical bases}

\begin{definition}\label {weak cb def}\em
Let $\Gamma_a$ be a partial type over a tuple $a$ (not necessarily finite). We say that a
hypeimaginary $e$ is a \em weak canonical base for $\Gamma_a$ \em if the following hold.\\
1) $e\in dcl(a)$ and $\Gamma_a$ does not fork over $e$.\\
2) $e\in bdd(e')$ whenever $e'\in bdd(a)$ is a hyperimaginary such that $\Gamma_a$ does not fork
over $e'$.
\end{definition}

\begin{example}\em
Let $L=\{E\}$ be a language for a 2-place relation $E$. Let $T$ be the complete $L$-theory saying
that $E$ is an equivalence relation with infinitely many equivalence classes each of which is
infinite. Let $\Gamma(x)\equiv E(x,a)\vee E(x,a')$ for some $a,a'\in\CC$ such that $\neg E(a,a')$.
Then $\Gamma(x)$ does not have a weak canonical base. To see this, assume otherwise. First, note
that clearly $\Gamma(x)$ doesn't fork over each of $a$ and $a'$. Thus by our assumption and the
definition of a weak canonical base, $e\in bdd(a)\cap bdd(a')$. So, we get a contradiction to the
fact that $\nonforkempty{a}{a'}$ and the fact that $\Gamma(x)$ forks over $\emptyset$.
\end{example}

We start by introducing the \em special partial types \em. First, we will say that a relation
$R(x,x')$ is \em generically transitive on a partial type $\pi(x)$ \em  if for all $
a',a,a''\models \pi$ if $\nonfork{a'}{a''}{a}$ and $R(a',a)\wedge R(a,a'')$, then $R(a',a'')$.

\begin{lemma}\label{weak Cb lemma}
Let $q(y,z),r(z,x)\in S(\emptyset)$ be such that $\exists x y z (q(y,z)\wedge r(z,x))$. Assume
$q(y,z)\vdash z\in acl(y)$. Let $p(x)=\exists z\ r(z,x)$ and let $a\models p$. Let $\Gamma_a$ be
defined by
$$\Gamma_a(y)\equiv\exists z\Bigg( q(y,z)\wedge r(z,a)\wedge \nonfork {y}{a}{z}\Bigg) .$$
\noindent Then\\
\noindent 1) $\Gamma_a(x)$ is a partial type.\\
\noindent 2) Let $\bar R_\Gamma$ be the relation defined by $\bar R_\Gamma(c,a,a')$ iff $c\models
\Gamma_{a}\wedge \Gamma_{a'}$ and $\nonfork{c}{a'}{a}$ and $\nonfork{c}{a}{a'}$. Then $\bar
R_\Gamma$ is type-definable and thus so is the relation $R_\Gamma(x,x')\equiv\exists y \bar
R_\Gamma(y,x,x')$ on $p^\CC$.

\noindent 3) If $q(y,z)\vdash z\in dcl(y)$ and $tp(d/a)$ is an amalgamation base for $(d,a)\models
r$, then $R_\Gamma$ is generically transitive on $p^\CC$.
\end{lemma}

\proof 1) is easy since in the definition of $\Gamma_a$ the complete type of $(z,a)$ is fixed. For
2) note that since $q(y,z)\vdash z\in acl(y)$, an easy forking computation shows that for all
$c,a,a'$ we have $\bar R_\Gamma(c,a,a')$ iff there exist $d,d'$ such that $q(c,d)\wedge
q(c,d')\wedge r(d,a)\wedge r(d',a')$ and
$\nonfork{c}{aa'}{d}\wedge\nonfork{c}{aa'}{d'}\wedge\nonfork{d}{a'}{a}\wedge \nonfork{d'}{a}{a'}$.
Again, since $q$ and $r$ are complete we get that $\bar R_\Gamma$ is type-definable. To prove 3),
assume $q(y,z)\vdash z\in dcl(y)$. Let $a,a',a''$ be such that $\nonfork{a'}{a''}{a}$ and assume
$R_\Gamma(a',a)$ and $R_\Gamma(a,a'')$ (so, clearly $a,a',a''\models p$). Now, by 2) we know that
$\bar R_\Gamma(x,a',a)$ is a partial type and clearly by its definition doesn't fork over $a$.
Likewise, the partial type $\bar R_\Gamma(x,a,a'')$ doesn't fork over $a$. Now, it will be
sufficient to show the following.

\begin{claim}\label{claim_1}

There are $c'\models\bar R_\Gamma(x,a',a)$ and $c''\models\bar R_\Gamma(x,a,a'')$ such that
$Lstp(c'/a)=Lstp(c''/a)$.

\end{claim}

\noindent This is sufficient since $\nonfork{c'}{a'}{a}$ for all $c'\models\bar R_\Gamma(x,a',a)$
and $\nonfork{c''}{a''}{a}$ for all $c''\models\bar R_\Gamma(x,a,a'')$, thus by the independence
theorem this will imply there exists $c^*\models \bar R_\Gamma(x,a',a)\wedge \bar
R_\Gamma(x,a,a'')$ with $\nonfork{c^*}{aa'a''}{a}$. In particular, $\nonfork{c^*}{aa'a''}{aa'}$,
and by the definition of $\bar R_\Gamma$, $\nonfork{c^*}{aa'a''}{a'}$. Likewise,
$\nonfork{c^*}{aa'a''}{a''}$. Hence $\bar R_\Gamma(a',a'')$.

\noindent\textbf{Proof of Claim \ref{claim_1}}. Since $q(y,z)\vdash z\in dcl(y)$, by the
observation on $\bar R_\Gamma$ in the proof of 2), there are $d'',c''$ such that $q(c'',d'')$,
$r(d'',a)$ and $r(d'',a'')$ and $\nonfork{d''}{a''}{a}$ and $\nonfork{d''}{a}{a''}$ and
$\nonfork{c''}{aa''}{d''}$. Likewise, there is $d'$ with $\nonfork{d'}{a'}{a}$ and
$\nonfork{d'}{a}{a'}$ and $r(d',a)$ and $r(d',a')$. Now, since $Lstp(d'/a)=Lstp(d''/a)$, there
exists $c'$ such that $tp(c'd'/a)=tp(c''d''/a)$ and $Lstp(c'/a)=Lstp(c''/a)$ and
$\nonfork{c'}{a'}{d'a}$. By the choice of $c'$, $\nonfork{c'}{a}{d'}$, hence by transitivity
$\nonfork{c'}{aa'}{d'}$. We conclude that $\bar R_\Gamma(c'',a,a'')$, $\bar R_\Gamma(c',a',a)$ and
$Lstp(c'/a)=Lstp(c''/a)$, as required.

\begin{definition}\label{special_def}\em
A partial type $\Gamma_a(y)$ as defined in Lemma \ref{weak Cb lemma} for some $q(y,z),r(z,x)\in
S(\emptyset)$ \em is called a special partial type \em if $q(y,z)\vdash z\in dcl(y)$ and $tp(d/a)$
is an amalgamation base for $(d,a)\models r$.
\end{definition}

\begin{remark}\em
Note that, in general, the class of special partial types properly contains the class of
amalgamation bases. To see this inclusion, let $r(z,x)$ be any complete type over $\emptyset$ such
that $tp(d/a)$ is an amalgamation base for $(d,a)\models r$, and let $q(y,z)=(y=z)\wedge(\exists x\
r(z,x))$. Then, if we apply the definition of $\Gamma_a$ in Lemma \ref{weak Cb lemma}, we get
$\Gamma_a(y)=tp(d/a)$ for $(d,a)\models r$. To justify properness, we give an example of a special
partial type which is not complete. Let $L=\{R\}$ and let $T$ be the $L$-theory of the random
graph. Let $a,b,c\in \CC$ be any three distinct elements such that say $R(b,c)$ and $R(c,a)$. Then
$\nonfork{bc}{a}{c}$. Let $q=tp(bc,c)$, $r=tp(c,a)$. Let
$$\Gamma_a(y_0y_1)\equiv\exists z\Bigg (q(y_0y_1,z)\wedge r(z,a)\wedge \nonfork {y_0y_1}{a}{z}\Bigg).$$ Then
clearly $\Gamma_a$ is a special partial type and $\Gamma_a(y_0y_1)$ is equivalent to $(y_0\neq
y_1)\wedge (y_0\neq a) \wedge (y_1\neq a)\wedge R(y_0,y_1)\wedge R(y_1,a)$. In particular,
$\Gamma_a$ is not complete.
\end{remark}

For proving the theorem we need the following well known fact (used for the construction of the
usual canonical base).

\begin{fact}\label{Cb fact} $[W, Lemma\ 3.3.1]$
Let $\pi(x)$ be a partial type over $\emptyset$ and let $R(x,x')$ be a type-definable relation over
$\emptyset$ that is reflexive, symmetric and generically transitive on $\pi^\CC$. Let $E_R$ be the
transitive closure of $R$ on $\pi^\CC$. Then $E_R$ is type-definable and for all $a,a'\models\pi$
we have $E_R(a,a')$ iff there exists $b\models\pi$ such that $R(a,b)$ and $R(b,a')$ and
$\nonfork{a}{b}{a'}$ and $\nonfork{b}{a'}{a}$.
\end{fact}

\begin{theorem}\label{wcb thm}
Let $\Gamma_a$ be a special partial type. Then $\Gamma_a$ has a weak canonical base.
\end{theorem}

\proof By Lemma \ref{weak Cb lemma} and Fact \ref{Cb fact} we know that the transitive closure of
$R_\Gamma$ (as defined in Lemma \ref{weak Cb lemma}), which we denote by $E_\Gamma$, is
type-definable and for all $a'\models tp(a)$ we have $E_\Gamma(a,a')$ iff there exists $b\models
tp(a)$ such that $R_\Gamma(a,b)$ and $R_\Gamma(b,a')$ and $\nonfork{a}{b}{a'}$ and
$\nonfork{b}{a'}{a}$. Let $e=a_{E_\Gamma}$. Clearly $e\in dcl(a)$. First we show $\Gamma_a$ doesn't
fork over $e$. Pick $a'\models tp(a)$ such that $\nonfork{a}{a'}{a_{E_\Gamma}}$ and $tp
(a'/a_{E_\Gamma})=tp(a/a_{E_\Gamma})$. In particular, $E_\Gamma(a,a')$. Let $b$ be as above. Then
easily $\nonfork{a}{b}{a_{E_\Gamma}}$ and thus if $c\models\Gamma_a\wedge \Gamma_b$ with
$\nonfork{a}{c}{b}$, then since $a_{E_\Gamma}\in dcl(b)$ we conclude that
$\nonfork{a}{c}{a_{E_\Gamma}}$. To prove 2) of \ref{weak cb def}, assume $\Gamma_a(x)$ doesn't fork
over some $e'\in bdd(a)$. Let $\sigma\in Aut(\CC/bdd(e'))$ and let $a'=\sigma(a)$. Pick $a^*$ such
that $tp(a^*/bdd(e'))=tp(a/bdd(e'))$ and such that $\nonfork{a^*}{aa'}{e'}$. By the independence
theorem, both $\Gamma_a\wedge\Gamma_{a^*}$ and $\Gamma_{a^*}\wedge \Gamma_{a'}$ doesn't fork over
$e'$. Since $e'\in bdd(a^*)\cap bdd(a)\cap bdd(a')$, $E_\Gamma(a,a')$. Thus $e\in bdd(e')$.\\

\begin{remark}\em
Definition \ref{special_def} of special partial types can be applied in the more general context of
hyperimaginaries. It is not hard to check that Theorem \ref{wcb thm} remains true in this context;
the main properties we need for that are the following. First, for the proof of Fact \ref{Cb fact}
we only need two properties besides standard forking computations; the first one is that
$\nonfork{a}{b}{c}$ if and only if for every $\phi=\phi(x,y)\in L$ and $k<\omega$ we have
$D(tp(a/c),\phi,k)=D(tp(a/bc),\phi,k)$ and the second property is the type-definability of the
$D(-,\phi,k)$-rank in the following sense: for every, possibly infinite, tuples of sorts $s_0,s_1$,
$\phi=\phi(x,y)\in L$ and $k,n<\omega$ the set $\{(a_0,a_1)\in \CC^{s_0}\times \CC^{s_1}\vert\
D(tp(a_1/a_0),\phi,k)\geq n\}$ is type-definable. These properties remains true in the
hyperimaginary context and thus so is Fact \ref{Cb fact}. For the proofs in this paper (and even
for knowing that special partial types are in fact types) we only need, in addition, the following
property: if $b_0,c_0$ are hyperimaginaries then for any fixed hyperimaginary sort $S_E$ (where $E$
is a type-definable equivalence relation over $\emptyset$), the set $\{(a,b,c) \vert\ a\in S^\CC_E\
, \nonfork{a}{b}{c}, \ tp(b,c)=tp(b_0,c_0)\}$ is type-definable (i.e. the union of the classes of
members of this set is type-definable).

\end{remark}

Here is a corollary of our main theorem. For simplicity we assume that $T$ is hypersimple (rather
than just simple). In the following, when we write $Cb(a/b)$, we mean the usual canonical base of
(the amalgamation base) $tp(a/bdd(b))$ (where $bdd(b)$ denotes the set of hyperimaginaries of
countable length whose type over $b$ is bounded). The assumption that $T$ is hypersimple implies
that such a canonical base exists as a set of imaginary elements and a type over an algebraically
closed set in $\CC^{eq}$ is an amalgamation base (since $bdd(A)$ is interdefinable with
$acl^{eq}(A)$ for every set $A\subseteq \CC^{eq}$).

\begin{corollary}\label {cor1}
Let $T$ be simple theory with elimination of hyperimaginaries and work in $\CC^{eq}$. Let $d,a$ be
some tuples (possibly infinite) and let $p\in S(d)$ be such that for $c\models p$, $d\in dcl(c)$.
Let
$$S=\{c \in p^\CC \vert \nonfork{c}{a}{d}\}.$$ Then there exists $c^*\in S$ such that
$\bigcap_{c\in S} acl(Cb(c/a))=acl(Cb(c^*/a))$.
\end{corollary}

\proof Let $\tilde a=acl(a)$. Let  $\Gamma_{\tilde a}$ be the special partial type over $a$ defined
by the types $r=tp(d,\tilde a)$ and $q=tp(c,d)$ for some $c\models p$. Note that $\{tp(c/\tilde a)
\vert\ c\in S\}=\{tp(c/\tilde a) \vert\ c\models \Gamma_{\tilde a}\}$, and clearly $S\subseteq
\Gamma_{\tilde a}^{\CC}$. By Theorem \ref{wcb thm}, there is a weak canonical base of
$\Gamma_{\tilde a}$, call it $e$. Let $c'\models \Gamma_{\tilde a}$ be such that
$\nonfork{c'}{a}{e}$, and let $e^*=Cb(c'/a)$. Then by the definition of the usual canonical base,
$e^*\in bdd(e)$. By the observation above, there exists $c^*\in S$ such that $e^*=Cb(c^*/a)$. Now,
to finish the proof it will be sufficient to show that $e\in bdd(Cb(c/a))$ for every $c\in S$.
Indeed, let $c\in S$, then $c\models \Gamma_{\tilde a}$. Let $e_c=Cb(c/a)$. Then
$\nonfork{c}{a}{e_c}$, and since $e$ is a
weak canonical base of $\Gamma_{\tilde a}$, we conclude $e\in bdd(e_c)$.\\


\begin{thebibliography}{www}
\bibitem[HKP]{www} B.Hart, B.Kim and A.Pillay, Coordinatization and canonical bases in simple theories,
Journal of Symbolic Logic, 65 (2000), pgs 293-309.
\bibitem[K]{www} B.Kim, Forking in simple unstable theories, Journal
of London Math. Society, 57 (1998), pgs 257-267.
\bibitem[KP]{www} B.Kim and A.Pillay, Simple theories, Annals of Pure and Applied Logic, 88, 1997 pgs 149-164.
\bibitem[W]{W} Frank~O. Wagner, Simple Theories, Academic Publishers, Dordrecht, The Netherlands, 2000.
\bibitem [S1]{www} Z.Shami, Countable hypersimple unidimensional theories,
J. London Math. Soc.  Volume 83, Issue 2 (2011), pgs. 309-332.

\end{thebibliography}
\end{document}